\documentclass[12pt]{amsart} 
\usepackage{amsmath,amssymb}      
\usepackage{ascmac}
\usepackage{mathrsfs}
\usepackage[abbrev]{amsrefs}
\usepackage{comment}

\usepackage{amsthm}
\usepackage{latexsym}
\usepackage{mathtools}
\usepackage{amsfonts}

\allowdisplaybreaks[1]
\newtheorem{thm}{Theorem}[section]

\newtheorem{prop}[thm]{Proposition}
\newtheorem{lem}[thm]{Lemma}

\newenvironment{pr}
   {{\noindent \bf Proof.   }}{\hfill \qed}

\def\<{\langle }

\newcommand{\eqsp}[1]{{\begin{equation}\begin{split}#1\end{split}\end{
equation}}}

\usepackage{graphicx,xcolor}

\begin{document}

% \maketitle

\begin{center}
\textbf{\Large{Higher-order derivative estimates for the heat equation on a smooth domain}}
\end{center}
\vskip10mm

\centerline{Yoshinori Furuto${}^*$ and Tsukasa Iwabuchi${}^{**}$}
\vskip5mm
\centerline{Mathematical Institute, Tohoku University}
\centerline{Sendai 980-8578 Japan}
\footnote[0]{E-mail: $^*$yoshinori.furuto.p3@dc.tohoku.ac.jp,
$^{**}$t-iwabuchi@tohoku.ac.jp}

\begin{center}
 \begin{minipage}{120mm}
\small \textbf{Abstract.} 
We consider the linear heat equation on a bounded domain and on an exterior domain. 
We study estimates of any order derivatives of the solution locally in time in the Lebesgue spaces.
We give a proof of the estimates in the end-point cases $p = 1, \infty$. 
We also obtain derivative estimates for the equation with the fractional Dirichlet Laplacian.
 \end{minipage}
\end{center}

% Section1 Introduction
\section{Introduction}

Let $d \geq 2$. Suppose that $\Omega$ is a bounded domain or an exterior domain of ${\mathbb R}^d$.
We consider the heat equation with the Dirichlet boundary condition.
\begin{equation} \label{eq:main}
\left \{
\begin{aligned}
\partial_t u - \Delta u &= 0
&& \quad \text{ in } \ (0, T) \times \Omega, 
\\
u &= 0
&& \quad \text{ on } \ (0, T) \times \partial \Omega, 
\\
u(0, \cdot) &= u_0(\cdot) 
&& \quad \text{ in } \ \Omega.
\end{aligned}
\right.
\end{equation}
Our aim is to obtain $L^p$ estimates of any order derivatives.

% Let us recall some classical results. 
% If $\Omega = \mathbb{R}^d$ and the boundary condition satisfies that $u(t, \cdot)$ converges to $0$ at infinity, then we can write a solution of~\eqref{eq:main} as
% \[
% \displaystyle u(t, x) = (G_t \ast u_0)(x) 
% = \int_{\mathbb{R}^d} \frac{1}{(4 \pi t)^{d/2}} \exp\left(-\frac{|(x - y)|^2}{4t} \right) u_0(y) \, \mathrm{d}y, 
% \]
% and it fulfills the following $L^p$-$L^p$ estimates
% \[
% \| \nabla^\ell u(t) \|_{L^p} \leq C t^{-\ell/2} \| u_0 \|_{L^p} \quad \text{ for any } \ell \in \mathbb{Z}.
% \]

If $\Omega$ is a sufficiently smooth bounded domain and $p \in (1, \infty)$, then there exists the heat semigroup $\{S_p(t)\}_{t \geq 0}$ on $L^p(\Omega)$ (see e.g.~\cite{book:Lunardi}). The solution $u(t) := S_p(t) u_0$, defined by the semigroup, satisfies the following $L^p$-$L^p$ estimates: 
There exists $C > 0$ such that
\begin{equation*}
\begin{aligned}
    \| \nabla^\ell u(t) \|_{L^p} &\leq C t^{-\ell / 2} \| u_0 \|_{L^p}
    \qquad \text{ for }  \ell \in \{ 0, 1, 2 \}, t > 0
\end{aligned}
\end{equation*}
where $\| \cdot \|_{L^p}$ is the norm of $L^p(\Omega)$. This is obtained by the elliptic estimates.
We notice that $u(t)$ decays exponentially in $L^p(\Omega)$ for large $t$, since the infimum of the spectrum of the Dirichlet Laplacian $-\Delta_{\mathrm D}$ is positive. 
Our aim is to consider higher-order derivative estimates, including $p = 1, \infty$, locally in time.

\begin{comment}
\begin{equation} \label{20240131-1}
    \begin{aligned}
        \| u(t) \|_p &\leq C \| u_0 \|_p, 
        \\
        \| \nabla u(t) \|_p &\leq C t^{-1/2} \| u_0 \|_p, 
        \qquad \text{ for } t > 0.
        \\
        \| \nabla^2 u(t) \|_p &\leq C t^{-1} \| u_0 \|_p,
    \end{aligned}
\end{equation}
In this estimate, we consider in local time.
The point of this estimate is to consider it in local time.
\end{comment}

In what follows, we suppose $u_0 \in C_c ^\infty (\Omega)$, which implies 
$S_2 (t) u_0 = S_p(t) u_0$ for any $p \in (1,\infty)$. 
We then abbreviate $p$, 
and often use the following notation  
\[
S(t)u_0 = S_p(t)u_0. 
\]

\vspace{3mm}

\begin{thm} \label{thm:p}
Let $d \geq 2, \gamma \in (\mathbb{Z}_{\geq 0})^d$. Suppose that $\Omega$ is a bounded domain or an exterior domain of ${\mathbb R}^d$ with $C^{|\gamma| + 2}$ boundary and $1 \leq p \leq \infty$. 
Then, a positive constant $C$ exists such that for every $u_0 \in C_c^\infty (\Omega)$, $u(t) := S(t) u_0$ satisfies that
\begin{equation} \label{eq:20240405-1}
\begin{aligned} 
    \|  \partial_x^\gamma u(t) \|_{L^p} &\leq C t^{-|\gamma| / 2} \| u_0 \|_{L^p}
    \quad \textup{ for all } t \in (0, 1), 
\end{aligned}
\end{equation}
where $\gamma = (\gamma_1, \gamma_2, \ldots, \gamma_d) \in (\mathbb{Z}_{\geq 0})^d$ and $|\gamma| = \gamma_1 + \gamma_2 + \cdots + \gamma_d$.
\end{thm}

\noindent 
{\bf Remark. } 
We provide a supplementary remark regarding the smoothness of the boundary in the assumption, which is two degree higher than the derivative order $|\gamma|$.
In the case when $p=1$, it is sufficient to impose $C^{|\gamma|}$ smoothness on the boundary, for the same reasons as the elliptic estimate in~\cite{book:Gilbarg-Trudinger}.
On the other hand, when $p = \infty$, we consider Sobolev space $W^{|\gamma|+1, r} (\subset W^{|\gamma|, \infty})$ in which the elliptic estimates hold.
Furthermore, we use the boundedness of $W^{|\gamma|+2, r}$ to claim the strong convergence of the $|\gamma|+1$-th order derivative, hence we impose the additional smoothness. 
% additional smoothness is required since we use $|\gamma|+1$-th derivative estimate in Sobolev space in the proof.

\vspace{3mm}

Consider the following generalization of the Theorem~\ref{thm:p} for fractional Dirichlet Laplacian. 
We define an operator $A^{\alpha / 2}$ and semigroup $e^{-tA^{\alpha / 2}}$ via the spectral theorem (see~\cite{paper:Iw-2018}).
\begin{equation*}
\begin{aligned}
&\mathcal{D} (A^{\alpha / 2}) = \left\{ f \in L^2(\Omega) \, \middle\vert \, \int^{\infty}_{-\infty} | \lambda |^\alpha \, \mathrm{d}\langle E_A(\lambda)f , f \rangle_{L^2(\Omega)} < \infty \right\}, &\quad t > 0, 
\\
&A^{\alpha / 2} = \int_{-\infty}^{\infty} \lambda^{\alpha / 2} \, \mathrm{d}E_A(\lambda), &\quad t > 0, 
\\
&e^{-tA^{\alpha / 2}} = \int_{-\infty}^{\infty} e^{-t \lambda^{\alpha / 2}}\, \mathrm{d}E_A(\lambda), &\quad t > 0, 
\end{aligned}
\end{equation*}
where $A$ is the a Dirichlet Laplacian on $L^2(\Omega)$ and $\{E_A(\lambda)\}_{\lambda \in \mathbb{R}}$ is a resolution of identity of $A$. 
Here note that the function $u(t,\cdot) = e^{-tA^{\alpha / 2}} u_0$ is a solution of the following equations.
\begin{equation*}
\left \{
\begin{aligned}
&\partial_t u + A^{\alpha/2} u = 0
&& \quad \text{ in } \ (0, T) \times \Omega, 
\\
&A^{n-1}u = \cdots = Au = u =  0
&& \quad \text{ on } \ (0, T) \times \partial \Omega, 
\\
&u(0, \cdot) = u_0(\cdot) 
&& \quad \text{ in } \ \Omega, 
\end{aligned}
\right.
\end{equation*}
where $n$ is a maximum integer satisfying $n \leq \alpha / 2$ (see section~3 in~\cite{book:Cazenave-Haraux} for the uniqueness of the solution).
In this case, the following theorem holds as higher-order derivative estimates for $e^{-tA^{\alpha / 2}}$ . 
\begin{thm} \label{thm:20250127-1}
Assume the same assumptions as in Theorem~\ref{thm:p}. Let $\alpha > 0$. 
Then, there exists a constant $C > 0$ such that for any $u_0 \in C_c^\infty(\Omega)$, it holds
\begin{equation} \label{eq:20250127-2}
\begin{aligned} 
    \| \partial_x^\gamma e^{-tA^{\alpha / 2}} u_0 \|_{L^p} &\leq C t^{-|\gamma| / \alpha} \| u_0 \|_{L^p}
    \quad \textup{ for all } t \in (0, 1).
\end{aligned}
\end{equation}
\end{thm}

\vspace{3mm}
% We will prove it by modifying the argument in~\cite{paper:IwMaTa-2018}.

%As known results about derivative estimates, 

\begin{comment}
書くか迷っているコメント
- u(t) = e^{-tA^{\alpha / 2}} u_0 が満たす方程式
- 積分核の各点評価により同様の微分評価が得られることが指摘されている
  →今回の手法は積分核の各点評価を経由せずに示しており, 別証明になっています
\end{comment}

\begin{comment}
To the best of our knowledge, 
Omega = R^d
Omega : bounded domain
derivative の順に書く場合
- u : 
p = [1, inf] でほぼ自明
- first derivative : 
p = [1, inf) は paper:IwTaMa-2021 の Section2 最後を引用
p = inf は Abe-Giga 
- second derivative : 
p = (1, inf) は Lunardi & G-T の手法が使える
p = inf は Abe-Giga
\end{comment}

\begin{comment}
We refer to some existing literature on derivative estimates.
When $\Omega = \mathbb R^d$, the inequalities~\eqref{eq:20240405-1} can be shown by a direct calculation (e.g.~\cite{book:Giga-Giga-Saal}). 
In the case when $\Omega$ is a smooth bounded domain and $p \in (1, \infty)$, $L^p$ estimates~\eqref{eq:20240405-1} can be obtained by using $L^p$ elliptic estimate and basic property of analytic semigroup (see~\cite{book:Lunardi} and~\cite{book:Gilbarg-Trudinger}).
In addition, \cite{paper:Iw-2018} shows a fist-order derivative estimate for $p = 1$.  
\end{comment}

We refer to some existing literature on derivative estimates.
When $\Omega = \mathbb{R}^d$, the inequalities~\eqref{eq:20240405-1} for $\gamma \in (\mathbb{Z}_{\geq 0})^d$ can be shown by a direct calculation.
In the case when $\Omega$ is a smooth bounded domain, the estimate of the solution itself without the derivatives is obvious from the property of the semigroup.
For $|\gamma| = 1, 2$, if $p \in (1, \infty)$, then~\eqref{eq:20240405-1} can be obtained by using $L^p$ elliptic estimate and basic property of analytic semigroup (see~\cite{book:Gilbarg-Trudinger} and~\cite{book:Lunardi}). In addition, \eqref{eq:20240405-1} for $|\gamma| \leq 2$ holds for $p = 1, \infty$ (see~\cite{paper:AbGi-2013} and~\cite{paper:FuIwKo-2024}).

A classical method for obtaining derivative estimates is pointwise estimates of integral kernels.
In the case of the heat kernel on a bounded domain, we refer {\`E}{\u{\i}}del'man and Ivasishen~\cite{paper:EI-1970} as a reference for a pointwise estimate of first-order derivatives.
However, the present approach, using resolvent estimates and elliptic estimates, derives $L^p$ estimates without relying on pointwise estimates.
This avoidance of pointwise estimates also applies to the case of fractional Dirichlet Laplacian.

% For the first derivatives, \eqref{20240409-1-b} holds for $p \in [1, \infty]$ (see~\cite{paper:Iw-2018},~\cite{paper:Shen-2005} and~\cite{paper:AbGi-2013}).
% For the second derivatives, if $p \in (1, \infty)$, then~\eqref{20240409-1-c} can be obtained by using $L^p$ elliptic estimate and basic property of analytic semigroup (see~\cite{book:Lunardi} and~\cite{book:Gilbarg-Trudinger}). In addition,~\eqref{20240409-1-c} holds for $p = \infty$ (see~\cite{paper:AbGi-2013}).
We give a few comments on our proof.
We prove only the case when $p = 1, \infty$.
The boundedness of without derivatives is due to the basic property of the $C_0$ semigroup.
For $p = 1$, we modify the method in~\cite{paper:IwMaTa-2018} to the higher-order resolvent of Laplacian to obtain a higher-order, in particular higher than second-order, derivative estimates.
As for $p = \infty$, we argue by contradiction used in~\cite{paper:AbGi-2013} for the Stokes equation. 
We extend an argument for higher order derivatives, introducing the following quantity
\[
N_\ell[u](t, x) := \sum_{s = 0}^{\ell} t^{s/2} |\nabla^s u(t, x)|.
\]
% and get the necessary boundedness by applying derivative estimates in a nonhomegeneous elliptic estimates repeatedly.
% Although not used in this paper, there is a possibility to show this theorem by direct estimate of the heat kernel using the method in \cite{paper:IsTa-2022}.

This paper is organized as follows.
In Section~2, we prepare a lemma and a proposition to prove the theorem for $p = \infty$.
In Sections~3 and~4, we prove the theorem for $p = 1, \infty$. 
In Section 5, we prove theorem~\ref{thm:20250127-1}.
In the Appendix, we show some lemmas used in the proof for $p = 1$.

\vskip3mm

\noindent {\bf Notation.} 
We denote by $A$ the Dirichlet Laplacian, which is defined by 
\[
\begin{cases}
D(A) = \{ f \in H^1_0 (\Omega) \, | \, \Delta f \in L^2 (\Omega)\}, 
\\
\displaystyle Af = - \sum _{ j=1}^d \partial _{x_j}^2 f , \quad f \in D(A) . 
\end{cases}
\]
We notice that $A$ is a self-adjoint operator on $L^2(\Omega)$ and by the argument in the paper ~\cite{paper:IwMaTa-2018}, $e^{-tA}$ can be defined based on the spectral multiplier. 
If $u_0 \in C_c^\infty (\Omega)$, then $S(t) u_0 = e^{-tA}u_0$, which is used in the proof for $p = 1$.
We denote by $\nabla^\ell$ the $\ell$th-order derivatives in all directions for any $\ell \in \mathbb{Z}_{\geq 0}$. That is, \eqref{eq:20240405-1} for all $\gamma$ with $|\gamma| = \ell$ is equivalent to
\[
\| \nabla^{\ell} u(t) \|_{L^p} \leq C t^{-\ell / 2} \| u_0 \|_{L^p}.
\]

% Section2 Preparation
\section{Preliminaries}

We prepare a lemma and a proposition for the proof when $p = \infty$.
We consider the estimate for the following quantity.
\begin{equation} \label{def:20241216-1}
N_\ell[u](t, x) := \sum_{s = 0}^{\ell} t^{s/2} |\nabla^s u(t, x)|.
\end{equation}

The first lemma considers an estimate of $N_\ell[u]$, and that for the Stokes equations for $\ell = 2$ is known~(see~\cite{paper:AbGi-2013}).
\begin{lem} \label{lemma:C2}
    Suppose that $\Omega$ is a bounded domain with $C^{\ell+1}$ boundary. Assume that $u_0 \in C_c^\infty(\Omega)$, and define $u(t) := S(t) u_0$.
    Then, % $u(t, \cdot) \in C^\ell(\Omega)$ for all $t > 0$ and
    \[
    \sup_{0 < t < 1}{\| N_\ell[u](t) \|_{L^\infty} }< \infty.
    \]
\end{lem}

To prove Lemma~\ref{lemma:C2}, we introduce higher order derivative estimates in a nonhomegeneous elliptic equation.

\begin{prop} \label{high-est-thm}
Suppose $\ell \in \mathbb N \cup \{0\}$. Let $\Omega \subset {\mathbb R}^d$ be a bounded domain or an exterior domain with $C^{\ell+2}$ boundary and $p \in (1, \infty)$.
If $u \in W^{2, p}(\Omega)$ and $f \in W^{\ell, p}(\Omega)$ satisfy
\begin{equation}
\left \{
\begin{aligned}
\Delta u &= f
&& \quad \textup{ in } \Omega, 
\\
u &= 0
&& \quad \textup{ on } \partial \Omega,
\end{aligned}
\right.
\end{equation}
then $u \in W^{\ell+2, p}(\Omega)$ and there exists $C = C(n, \ell, p)$ such that
\[
\| u \|_{W^{\ell+2, p}(\Omega)} \leq C ( \| u \|_{L^p(\Omega)} + \| f \|_{W^{\ell, p}(\Omega)} ).
\]
\end{prop}

Proof idea can be found in the book by Gilbarg-Trudinger~\cite{book:Gilbarg-Trudinger} (Theorem 9.19), and the proof is similar to Theorem~A.1. in the paper~\cite{paper:FuIwKo-2024} if $\Omega$ is an exterior domain.

\begin{proof}[{\noindent \bf Proof of Lemma~\ref{lemma:C2}.}]
    Let $d < r < \infty$, and $-A_r$ be the generator of the heat semigroup in $L^r(\Omega)$.
    We will show that there exists $C = C(\Omega, r) > 0$ such that
    \begin{equation} \label{equation:LemmaGoal}
    \begin{split}
        \sup_{0 < t < 1}{ \sum_{s = 0}^{\ell} t^{s/2} \| \nabla^s u(t, x) \|_{W^{1, r}} }
        \leq C ( \| u_0 \|_{D(A_r)} + \| \Delta u_0 \|_{D(A_r)} 
          + \cdots + \| \Delta^{\lceil (\ell+1)/2 \rceil} u_0 \|_{D(A_r)} ),
    \end{split}
    \end{equation}
    where $\| u_0 \|_{D(A_r)} := \| u_0 \|_{L^r} + \| A_r u_0 \|_{L^r}$, $\lceil \alpha \rceil$ is the smallest integer greater than or equals to $\alpha$, for $\alpha \in \mathbb{R}$.

    By the property of the analytic semigroup $S_r(t) = e^{-t A_r}$, we have
    \[
    \begin{aligned}
        \sup_{0 < t < 1}{\| u(t) \|_{D(A_r)}} 
        \leq C \| u_0 \|_{D(A_r)},
    \end{aligned}
    \]
    where $C = C(\Omega, r)$. Thus we have proved
    % \[
    % \sup_{0 < t < 1}{ (\| u(t) \|_{W^{1, r}} + \| \nabla u(t) \|_{W^{1, r}} }
    % \leq C \| u_0 \|_{D(A_r)}
    % \]
    \[
    \sup_{0 < t < 1}{ ( \| u(t) \|_{L^r} + \| \nabla u(t) \|_{L^r} + \| \nabla^2 u(t) \|_{L^r} ) }
    \leq C \| u_0 \|_{D(A_r)}
    \]
    as the $D(A_r)$ norm and the $W^{2, r}(\Omega)$ norm are equivalent by the elliptic estimate.
    
    % The terms
    % \[
    % \| u(t) \|_{W^{1, r}}, \quad
    % t^{1/2} \| \nabla u(t) \|_{W^{1, r}}, \quad 
    % t \| \nabla^2 u(t) \|_{L^r}, \quad
    % t \| \partial_t u(t) \|_{W^{1, r}}
    % \]
    % are less than or equal to the left hand side above.
    %The left hand side controls the terms
    % \[
    % || u(t) ||_{W^{1, r}(\Omega)}, t^{1/2} ||a \nabla u(t) ||_{W^{1, r}(\Omega)}, 
    % t || \nabla^2 u(t) ||_{L^r(\Omega)}, t || \partial_t u(t) ||_{W^{1, r}(\Omega)}.
    % \]

    We use Proposition~\ref{high-est-thm} repeatedly to estimate higher order derivatives. Applying Proposition~\ref{high-est-thm} for $(u, f) = (u, \Delta u)$, we have
    \begin{equation*}
    \begin{aligned}
    \| u(t) \|_{W^{s, r}} 
    \leq C ( \| u(t) \|_{L^r} + \| \Delta u(t) \|_{W^{s-2, r}}), 
    \end{aligned}
    \end{equation*}
    and similarly, applying it for $(u, f) = (\Delta u, \Delta^2 u)$, then
    \begin{equation*}
    \begin{aligned}
    \| \Delta u(t) \|_{W^{s-2, r}} 
    = C ( \| \Delta u(t) \|_{L^r} + \| \Delta^2 u(t) \|_{W^{s-4, r}} ).
    \end{aligned}
    \end{equation*}
    Note that $\Delta u(t) \in D(A)$.
    Therefore by induction, for every $j \in \mathbb{N}$, in the case when $s = 2j$ , 
    \begin{equation*}
    \begin{aligned}
    \| u(t) \|_{W^{2j, r}}
    =&\, C ( \| u(t) \|_{L^r} + \| \Delta u(t) \|_{W^{2j-2, r}} )
    \\
    \leq&\, C ( \| u(t) \|_{L^r} + \| \Delta u(t) \|_{L^r} + \| \Delta^2 u(t) \|_{W^{2j-4, r}} )
    \\
    &\vdots
    \\
    \leq&\, C ( \| u(t) \|_{L^r} + \| \Delta u(t) \|_{L^r} + \cdots + \| \Delta^j u(t) \|_{L^r} ), 
    \end{aligned}
    \end{equation*}
    and when $s = 2j + 1$, 
    \begin{equation*}
    \begin{aligned}
    \| u(t) \|_{W^{2j+1, r}} 
    \leq&\, C ( \| u(t) \|_{L^r} + \| \Delta u(t) \|_{L^r} + \cdots + \| \Delta^j u(t) \|_{W^{1, r}} ).
    % \leq& C ( \| u(t) \|_{L^r} + \| \Delta u(t) \|_{L^r} + \cdots + \| \Delta^{s+1} u(t) \|_{L^r} ).
    \end{aligned}
    \end{equation*}
    
    In addition, due to the elliptic estimate and the property of the semigroup, it holds
    \[
    \| \Delta^j u(t) \|_{W^{2, r}}
    \leq C \| A^j e^{tA} u_0 \|_{D(A_r)}
    = C \| e^{tA} (A^j u_0) \|_{D(A_r)}
    \leq C \| \Delta^j u_0 \|_{D(A_r)}.
    \]
    We therefore obtain, regardless of whether s is even or odd,
    \[
    \| \nabla^s u(t) \|_{L^r}
    \leq C ( \| u_0 \|_{D(A_r)} + \| \Delta u_0 \|_{D(A_r)} 
      + \cdots + \| \Delta^{\lceil s/2 \rceil} u_0 \|_{D(A_r)} ).
    \]
    % \[
    % \| \nabla^s u(t) \|_{L^r}
    % \leq C \| \Delta^{\lceil (s-1)/2 \rceil} u_0 \|_{D(A_r)} .
    % \]
    By summing these inequalities from $s = 0$ to $s = \ell + 1$, we have
    \begin{equation*}
    \begin{aligned}
    \sup_{0 < t < 1}{ \sum_{s = 0}^{\ell} t^{s/2} \| \nabla^s u(t, x) \|_{W^{1, r}} }
    \leq&\, C \sum_{s = 0}^{\ell + 1} \| \nabla^s u(t, x) \|_{L^r}
    \\
    \leq&\, C ( \| u_0 \|_{D(A_r)} + \| \Delta u_0 \|_{D(A_r)} 
      + \cdots + \| \Delta^{\lceil (\ell+1)/2 \rceil} u_0 \|_{D(A_r)} ), 
    \end{aligned}
    \end{equation*}
    which proves~\eqref{equation:LemmaGoal}.
    Here the assumption that $\Omega$ is a $C^{\ell+1}$ domain is needed to apply Proposition~\ref{high-est-thm}.
\end{proof}

We introduce the following proposition on the uniqueness for the heat equation on the half space, which is used in the proof of Theorem~\ref{thm:p}.
The proof is similar to that for the whole space in~\cite{book:Giga-Giga-Saal}.

\begin{prop}[\cite{paper:FuIwKo-2024}] \label{prop:0126-1}
Let $\Omega = \mathbb R^d_+$. 
Suppose that $u$ satisfies
\[
\int_0^T \int_{\Omega} u (\partial_t \phi + \Delta \phi) {\mathrm d}x {\mathrm d}t 
= 0
\]
for any
$
\phi \in C([0, T] ; W^{1,1}_0(\Omega) \cap W^{2, 1}(\Omega) \cap C(\Omega)) 
\cap C^1([0, T] ; L^1(\Omega))
$
with $| \nabla \phi(x)| \leq C(1 + |x|^2)^{-d/2}$.
Then $u$ is $0$ almost everywhere.
\\
\end{prop}
% We provide a proof of Proposition~\ref{prop:0126-1} in Appendix B.

\section{$L^1$ estimate in Theorem~\ref{thm:p}} 

We follow the argument as in the paper~\cite{paper:JeNa-1994} and~\cite{paper:JeNa-1995} (see also~\cite{paper:IwMaTa-2018}). 
Instead of $S(t)u_0$, we consider more general problem, which is the boundedness of the 
spectral multiplier $\varphi(t A)$, 
where $\varphi : \mathbb R \to \mathbb R$ 
satisfies ${\rm supp \, }\varphi \subset [-1,\infty)$. 
We will show 
\[
\| \nabla^\ell \varphi(tA) f \|_{L^1}
\leq C t^{-\ell / 2} \| f \|_{L^p}, 
\quad t \in (0, 1)
\]
for any $f \in C_c^\infty (\Omega)$, where $C > 0$ is independent of $t$ and $f$, and dependent on $\ell$ and $\varphi$.

Following the proof of Theorem~2.1 in \cite{paper:IwMaTa-2018}, 
we decompose $\Omega$ by using cubes whose side length is $t^{\frac{1}{2}}$, and apply the H\"older 
inequality to have that 
\begin{equation}\label{eq:0131-1}
\| \nabla ^\ell \varphi (tA) f \|_{L^1} 
\leq t^{\frac{d}{4}} 
\| \nabla ^\ell \varphi(tA) f \|_{\ell^1 (L^2)_t},
\end{equation}
where 
\[
\| \nabla ^\ell \varphi(tA) f \|_{\ell^1 (L^2)_t}
:= \sum_{n \in {\mathbb Z}^d} \| \nabla ^\ell \varphi(tA) f \|_{L^2 (C_t(n))}, 
\]
and $C_t(n)$ is a cube with the side length $t^{1/2}$ centered at $t^{1/2} n \in {\mathbb R}^d$.
Fix a sufficiently large positive integer $\beta$ and $M$. Introduce $\widetilde \varphi \in \mathcal S(\mathbb R), \psi \in C_c^\infty(\mathbb{R})$ such that  
\[
\begin{aligned}
& \widetilde\varphi(\lambda) \coloneqq (1 + \lambda)^\beta \varphi(\lambda), 
\\
&\psi(\mu) \coloneqq \rho(\mu) \mu^{-1} \widetilde{\varphi}\left( \mu^{-\frac{1}{M}} - 1 \right), 
\end{aligned}
\]
where $\rho$ satisfies
\begin{equation*}
\rho \in C_c^\infty (\mathbb{R}), 
\quad
\rho(\mu) = 
\left \{
\begin{aligned}
1 \quad &(0 \leq \mu \leq 2)
\\
0 \quad &(\mu \leq -1, 3 \leq \mu).
\end{aligned}
\right.
\end{equation*}
% \color{red}{次のように書ける理由を明記する} \color{black}
Then, we notice that $\widetilde{\varphi}$ and $\psi$ satisfies
\[
\nabla^\ell \widetilde{\varphi}(tA) 
= \nabla^\ell (1 + tA)^{-M} \psi\left( (1 + tA)^{-M} \right) 
=: \nabla^\ell R_t \psi(R_t).
\]

We can write 
\begin{equation}\notag 
\begin{split}
& \| \nabla ^\ell \varphi(tA) f \|_{\ell^1 (L^2)_t}
\\
=& \| \nabla ^\ell \widetilde \varphi(tA) 
(1+tA)^{-\beta} f \|_{\ell^1 (L^2)_t}
\\
\leq & C 
\Big(  \| \nabla ^\ell \widetilde \varphi(tA) \| _{L^2 \to L^2} 
+ t ^{-\frac{d}{4}} 
 ||| \nabla ^\ell \widetilde \varphi(tA) |||_\alpha 
  ^{\frac{d}{2\alpha}} 
 \| \nabla ^\ell \widetilde \varphi(tA) \| _{L^2 \to L^2} ^{1-\frac{d}{2\alpha}} 
\Big)
\| (1+tA)^{-\beta} f \|_{\ell^1 (L^2)_t}
\end{split}
\end{equation}
for any $\alpha > d/2$, where 
\[
||| \nabla ^\ell \widetilde \varphi(tA) |||_\alpha :=
\sup_{n \in {\mathbb Z}^d} \left\| | \cdot - t^{1/2} n |^\alpha \nabla ^\ell \widetilde \varphi(tA) \chi_{C_t(n)}  \right\|_{L^2 \to L^2} . 
\]

Consider estimates of the norms that appear in the last inequality. Initially, we can prove the following.
\begin{prop} \label{prop:20241205-1}
Let $M$ and $\ell$ be nonnegative integers satisfying $0 \leq \ell \leq 2M$. Define $R_{M, t} = (1 + tA)^{-M}$ and assume $\widetilde \varphi \in \mathcal{S}(\mathbb{R})$. Then, there exists $C > 0$ such that
\[
\begin{aligned}
\| \nabla ^\ell R_{M, t} f \|_{L^2} \leq C t^{-\ell/2} \| f \|_{L^2}, 
\\
\| \nabla ^\ell \widetilde \varphi(tA) f\| _{L^2} \leq C t^{-\ell / 2} \| f \|_{L^2}
\end{aligned}
\]
for all $f \in L^2(\Omega)$ and $t \in (0, 1)$.
\end{prop}
\begin{pr}
The elliptic estimate in $L^2(\Omega)$ and Theorem~6.1 in~\cite{paper:IwMaTa-2018} shows that for any even number $\ell$, 
\begin{equation*}
\begin{aligned}
\| \nabla^\ell R_{M, t} f\|_{L^2}
\leq&\, C \| \nabla^{\ell - 2} \Delta  R_{M, t} f \|_{L^2} 
  + C \| R_{M, t} f \|_{L^2}
\\
\leq&\, C \| \nabla^{\ell - 2} t^{-1}((1 + tA) - 1) R_{M, t} f \|_{L^2} 
  + C \| R_{M, t} f \|_{L^2}
\\
\leq&\, C t^{-1} \| \nabla^{\ell - 2} (1+tA) R_{M, t} f \|_{L^2} 
  + C t^{-1} \| \nabla^{\ell - 2} R_{M, t} f \|_{L^2}  + C \| R_{M, t} f \|_{L^2}
\\
\leq&\, C t^{-1} \| \nabla^{\ell - 2} R_{M-1, t} f \|_{L^2} 
  + C t^{-1} \| \nabla^{\ell - 2} R_{M, t} f \|_{L^2}  + C \| R_{M, t} f \|_{L^2}
\\
& \vdots
\\
\leq&\, C t^{-\ell / 2} \| R_{M - \ell/2, t} f \|_{L^2} 
  \\
  &\, + C (t^{-1} \| \nabla^{\ell - 2} R_{M, t} f \|_{L^2} + \cdots + t^{-\ell / 2} \| R_{M, t} f \|_{L^2})
  \\
  &\, + C (1 + t^{-1} + \cdots + t^{-\frac{\ell}{2}}) \| R_{M, t} f \|_{L^2}
\\
\leq&\, C t^{-\frac{\ell}{2}} \| f \|_{L^2}, 
\\
\| \nabla^\ell \widetilde \varphi(tA) f\| _{L^2} 
=&\, \| \nabla^\ell R_{M, t} \cdot (1 + tA)^M \widetilde \varphi(tA) f\|_{L^2}
\\
\leq&\, C \| \nabla^\ell R_{M, t} \|_{L^2 \to L^2} \| (1 + tA)^M \widetilde \varphi(tA) f\|_{L^2}
\\
\leq&\, C t^{-\frac{\ell}{2}} \| f \|_{L^2}.
\end{aligned}
\end{equation*}
% \[
% \begin{aligned}
% \| \nabla^\ell R_{M, t} f\|_{L^2}
% \leq& C \| \nabla^{\ell - 2} \Delta  R_{M, t} f \|_{L^2} 
%   + C \| R_{M, t} f \|_{L^2}
% \\
% \leq& C t^{-1} \| \nabla^{\ell - 2} (1 + tA) R_{M, t} f \|_{L^2} 
%   + C t^{-1} \| \nabla^{\ell - 2} R_{M, t} f \|_{L^2}  + C \| R_{M, t} f \|_{L^2}
% \\
% & \vdots
% \\
% \leq& C t^{-\ell / 2} \| (1 + tA)^{\ell / 2} R_{M, t} f \|_{L^2} 
%   + C (t^{-1} \| \nabla^{\ell - 2} R_{M, t} f \|_{L^2} + \cdots + t^{-\ell / 2} \| R_{M, t} f \|_{L^2})
%   \\
%   &\, + C (1 + t^{-1} + \cdots + t^{-\frac{\ell}{2} - 1}) \| R_{M, t} f \|_{L^2}
% \\
% \leq& C t^{-\ell / 2} \| f \|_{L^2}, 
% \\
% \| \nabla^\ell \widetilde \varphi(tA) f\| _{L^2} 
% =& \| \nabla^\ell R_{M, t} \cdot (1 + tA)^M \widetilde \varphi(tA) f\|_{L^2}
% \\
% \leq& C \| \nabla^\ell R_{M, t} \|_{L^2 \to L^2} \| (1 + tA)^M \widetilde \varphi(tA) f\|_{L^2}
% \\
% \leq& C t^{-\frac{\ell}{2}} \| f \|_{L^2}.
% \leq& C ( \| A \nabla^{\ell-2} \widetilde \varphi(tA) f\| _{L^2} + \| \nabla^{\ell-2} \widetilde \varphi(tA) f\| _{L^2} )
% \\
% &\vdots
% \\
% \leq& C (\| A^{\frac{\ell}{2}} \widetilde \varphi(tA) f\| _{L^2} + \| A^{\frac{\ell}{2} - 1} \widetilde \varphi(tA) f\| _{L^2} + \cdots + \| \widetilde \varphi(tA) f\| _{L^2} )
% \\
% \leq& C (t^{-\frac{\ell}{2}} \| (tA)^{\ell / 2} \widetilde \varphi(tA) f\| _{L^2} + t^{-\frac{\ell}{2} + 1} \| (tA)^{\frac{\ell}{2} - 1} \widetilde \varphi(tA) f\| _{L^2} + \cdots + \| \widetilde \varphi(tA) f\| _{L^2} )
% \\
% \leq& C t^{-\frac{\ell}{2}} \| f \|_{L^2}.
% \end{aligned}
% \]
We used the fact $(1 + \cdot)^M \widetilde \varphi(\cdot) \in \mathcal{S}(\mathbb{R})$ in the last inequality. 
Combined with Lemma~7.1 in~\cite{paper:IwMaTa-2018}, it can be proved that these inequalities also hold when $\ell$ is odd.
\end{pr}

\vspace{3mm}

Now back to estimate of
$ ||| \nabla ^\ell \widetilde \varphi(tA) |||_\alpha $.
For the same argument as in the proof of Lemma~7.1 in \cite{paper:IwMaTa-2018}, we aim to show that
\begin{equation} \label{eq:2024-12-03-01}
\left\| | \cdot - t^{1/2} n |^\alpha \nabla ^\ell R_{M, t} e^{i \tau R_{M, t}} \chi_{C_t(n)}  \right\|_{L^2 \to L^2} \leq C (1 + |\tau|^\alpha) t^{\frac{\alpha}{2} - \frac{\ell}{2}}
\end{equation}

We focus on the case when $\alpha = 1$. 
For an operator $L$ on $L^2(\Omega)$, define a commutator
\[
\begin{aligned}
\mathrm{Ad}^0 (L) &:=  L, 
\\
\mathrm{Ad}^k (L) &:= \mathrm{Ad}^{k-1} \left( (x_j - t^{1/2}n_j) L - L (x_j - t^{1/2}n_j) \right), 
\end{aligned}
\]
where $x_j, n_j$ denote the $j$-th components of $x \in \mathbb{R}^d, n \in \mathbb{Z}^d$.
Then, notice that 
\[
\mathrm{Ad}^k (tA) = 
\begin{cases}
tA & \text{ if } k = 0, \\
2t \partial_{x_j} & \text{ if } k = 1, \\
-2t & \text{ if } k = 2, \\
0 & \text{ if } k \geq 3. \\
\end{cases}
\]

For example, when $\ell = 2$, the expression can be expanded as follows by using commutators.
\[
\begin{aligned}
& (x_j - t^{1/2} n_j) \partial_{x_i} \partial_{x_{i'}} R_{M, t} e^{i \tau R_{M, t}}
\\
=&\, \delta_{i, j} \partial_{x_{i'}} R_{M, t} e^{i \tau R_{M, t}}
  + \partial_{x_i} (x_j - t^{1/2} n_j) \partial_{x_{i'}} R_{M, t} e^{i \tau R_{M, t}}
\\
=&\, \delta_{i, j} \partial_{x_{i'}} R_{M, t} e^{i \tau R_{M, t}} 
  + \partial_{x_i} \delta_{i'j} R_{M, t} e^{i \tau R_{M, t}}
  + \partial_{x_i} \partial_{x_{i'}} (x_j - t^{1/2} n_j)  R_{M, t} e^{i \tau R_{M, t}}
\\
=&\, \delta_{i, j} \partial_{x_{i'}} R_{M, t} e^{i \tau R_{M, t}} 
  + \partial_{x_i} \delta_{i'j} R_{M, t} e^{i \tau R_{M, t}}
  + \partial_{x_i} \partial_{x_{i'}} \mathrm{Ad}^1 (R_{M, t}) e^{i \tau R_{M, t}}
  \\
  &\, + \partial_{x_i} \partial_{x_{i'}} R_{M, t} (x_j - t^{1/2} n_j) e^{i \tau R_{M, t}}
\\ 
=&\, \delta_{i, j} \partial_{x_{i'}} R_{M, t} e^{i \tau R_{M, t}} 
  + \partial_{x_i} \delta_{i'j} R_{M, t} e^{i \tau R_{M, t}} 
  + \partial_{x_i} \partial_{x_{i'}} \mathrm{Ad}^1 (R_{M, t}) e^{i \tau R_{M, t}} 
  \\
  &\, + \partial_{x_i} \partial_{x_{i'}} R_{M, t} \mathrm{Ad}^1 (e^{i \tau R_{M, t}})
  + \partial_{x_i} \partial_{x_{i'}} R_{M, t} e^{i \tau R_{M, t}} (x_j - t^{1/2} n_j). 
\end{aligned}
\]
Therefore, for $\ell \geq 1$, $j = 1, \ldots, d$ and $f \in L^2(\Omega)$, we can estimate the left-hand side of~\eqref{eq:2024-12-03-01}. 
\[
\begin{aligned}
&\left\| | x_j - t^{1/2} n_j | \nabla ^\ell R_{M, t} e^{i \tau R_{M, t}} \chi_{C_t(n)} f \right\|_{L^2}
\\
\leq&\, \ell \| \nabla^{\ell - 1} R_{M, t} e^{i \tau R_{M, t}} \chi_{C_t(n)} f \|_{L^2} + 
  \| \nabla^\ell \mathrm{Ad}^1 (R_{M, t}) e^{i \tau R_{M, t}} \chi_{C_t(n)} f \|_{L^2} 
  \\
  &\, + \| \nabla^\ell R_{M, t} \mathrm{Ad}^1 (e^{i \tau R_{M, t}}) \chi_{C_t(n)} f \|_{L^2} + 
  \| \nabla^\ell R_{M, t} e^{i \tau R_{M, t}} (x_j - t^{1/2} n_j)  \chi_{C_t(n)} f \|_{L^2}
\\
=:& \, I + II + III + IV.
\end{aligned}
\]

For $II$, we see that $\mathrm{Ad}^1 (R_{M, t})$ can be decomposed as below.
\[
\begin{aligned}
\mathrm{Ad}^1 ( (1 + tA)^{-1} ) 
&= (1+tA)^{-1}(x_j - t^{1/2} n_j) - (x_j - t^{1/2} n_j) (1+tA)^{-1}
\\
&= (1+tA)^{-1} \left[ (x_j - t^{1/2} n_j) (1 + tA) -  (1 + tA) (x_j - t^{1/2} n_j) \right] (1+tA)^{-1}
\\
&= (1+tA)^{-1} (2t \partial_{x_j}) (1+tA)^{-1},
\\
\mathrm{Ad}^1 ( (1 + tA)^{-2} ) 
&= (1+tA)^{-2}(x_j - t^{1/2} n_j) - (x_j - t^{1/2} n_j) (1+tA)^{-2}
\\
&= (1+tA)^{-2} \left[ (x_j - t^{1/2} n_j) (1 + tA) -  (1 + tA) (x_j - t^{1/2} n_j) \right] (1+tA)^{-1}
  \\
  &\, + (1+tA)^{-1} \left[ (x_j - t^{1/2} n_j) (1 + tA) -  (1 + tA) (x_j - t^{1/2} n_j) \right] (1+tA)^{-2}
\\
&= (1+tA)^{-2} (2t \partial_{x_j}) (1+tA)^{-1} + (1+tA)^{-1} (2t \partial_{x_j}) (1+tA)^{-2}.
\end{aligned}
\]
Similarly, we can check
\[
\mathrm{Ad}^1 ( R_{M, t} )
= \mathrm{Ad}^1 ( (1 + tA)^{-M} )
= \sum_{m_1 + m_2 = M+1} (1+tA)^{-m_1} (2t \partial_{x_j}) (1+tA)^{-m_2}, 
\]
where $m_1$ and $m_2$ in the above sum satisfies $m_1, m_2 \geq 1$.

In addition, we can prove $\ell-$th derivatives estimate for $\ell \leq 2(m_1 + m_2) - 1$
\[
\| \nabla^\ell (1+tA)^{-m_1} (2t \partial_{x_j}) (1+tA)^{-m_2} f \|_{L^2}
\leq C t^{-\frac{\ell - 1}{2}} \| f \|_{L^2}, t \in (0, 1), 
\]
by using the following estimate and induction on $\ell$. 
In fact, the case when $\ell = 0, 1$ follows obviously. If we assume the inequalities up to $(\ell-1)$-th derivative order, then
\[
\begin{aligned}
 & \| \nabla^\ell (1+tA)^{-m_1} (2t \partial_{x_j}) (1+tA)^{-m_2} f \|_{L^2}
\\
=&\, \| \nabla^\ell A^{-1} A (1+tA)^{-m_1} (2t \partial_{x_j}) (1+tA)^{-m_2} f \|_{L^2}
\\
\leq&\, C \| \nabla^{\ell-2} A (1+tA)^{-m_1} (2t \partial_{x_j}) (1+tA)^{-m_2} f \|_{L^2}
\\
  &\, + C \| (1+tA)^{-m_1} (2t \partial_{x_j}) (1+tA)^{-m_2} f \|_{L^2}
\\
\leq&\, C \| \nabla^{\ell-2} t^{-1} \left( (1+tA)^{-m_1 + 1} - (1+tA)^{-m_1} \right) (2t \partial_{x_j}) (1+tA)^{-m_2} f \|_{L^2}
\\
  &\, + C t \cdot t^{-1/2} \| f \|_{L^2}
\\
\leq&\,  C t^{-1}  \| \nabla^{\ell-2} (1+tA)^{-m_1 + 1} (2t \partial_{x_j}) (1+tA)^{-m_2} f \|_{L^2}
\\
  &\, + C t^{-1}  \| \nabla^{\ell-2} (1+tA)^{-m_1} (2t \partial_{x_j}) (1+tA)^{-m_2} f \|_{L^2}
\\
  &\, + C t^{1/2} \| f \|_{L^2}
\end{aligned}
\]
Here, the second inequality above is derived from Proposition~\ref{high-est-thm}. 
This proves the case of $\ell$-th order. 
This derivative estimate implies 
\[
II \leq \| \nabla^\ell \mathrm{Ad}^1 (R_{M, t}) \|_{L^2 \to L^2}  \| e^{i \tau R_{M, t}} \|_{L^2 \to L^2}  \| \chi_{C_t(n)} f \|_{L^2} \leq C t^{-\frac{\ell-1}{2}} \| f \|_{L^2}, 
\quad t \in (0, 1).
\]

Similarly, Proposition~\ref{prop:20241205-1} and Appendix in the present paper and Lemma~7.1 in~\cite{paper:IwMaTa-2018} imply the estimate of  $I, III, IV$. For $t \in (0, 1)$, 
\begin{equation*}
\begin{aligned}
&I \leq C \| \nabla^{\ell - 1} R_{M, t} \|_{L^2 \to L^2}  \| e^{i \tau R_{M, t}} \|_{L^2 \to L^2}  \| \chi_{C_t(n)} f \|_{L^2} \leq C t^{-\frac{\ell-1}{2}} \| f \|_{L^2}, 
\\
&III \leq \| \nabla^\ell R_{M, t} \|_{L^2 \to L^2}  \| \mathrm{Ad}^1 (e^{i \tau R_{M, t}}) \|_{L^2 \to L^2}  \| \chi_{C_t(n)} f \|_{L^2} \leq C t^{-\frac{\ell-1}{2}} (1 + |\tau|) \| f \|_{L^2}, 
\\
&IV \leq \| \nabla^\ell R_{M, t} \|_{L^2 \to L^2}  \| e^{i \tau R_{M, t}} \|_{L^2 \to L^2}  \| (x_j - t^{1/2} n_j)  \chi_{C_t(n)} f \|_{L^2} \leq C t^{-\frac{\ell-1}{2}} \| f \|_{L^2}.
\end{aligned}
\end{equation*}

Therefore we obtain~\eqref{eq:2024-12-03-01} for $\alpha = 1$
\begin{equation*}
\begin{aligned}
\left\| | \cdot - t^{1/2} n |^1 \nabla ^\ell R_{M, t}  e^{i \tau R_{M, t}} \chi_{C_t(n)}  \right\|_{L^2 \to L^2}
\leq&\, C t^{\frac{1 - \ell}{2}} (1 + |\tau|), 
\\    
||| \nabla ^\ell \widetilde \varphi(tA) |||_1 
\leq&\, C t ^{\frac{1 - \ell}{2}}, 
\end{aligned} 
\quad 0 < t < 1.
\end{equation*}
for $t \in (0, 1)$. 
For general $\alpha \in \mathbb{N}$, we would consider $\mathrm{Ad}^k$ with $k \geq 1$, and we get~\eqref{eq:2024-12-03-01} by using the equality in Appendix~\ref{sec:20241217-4}.
\begin{equation*}
\begin{aligned}
\left\| | \cdot - t^{1/2} n |^\alpha \nabla ^\ell R_{M, t}  e^{i \tau R_{M, t}} \chi_{C_t(n)}  \right\|_{L^2 \to L^2}
\leq&\, C t^{\frac{\alpha - \ell}{2}} (1 + |\tau|^\alpha), 
\\
||| \nabla ^\ell \widetilde \varphi(tA) |||_\alpha
\leq&\, C t ^{\frac{\alpha - \ell}{2}}, 
\end{aligned} 
\quad 0 < t < 1.
\end{equation*}
% \[
% \begin{aligned}
% \mathrm{Ad}^k ( R_{M, t} )
% =& \sum_{i = 0}^{k-1} \sum_{m_1+m_2 = M+1} 
%   && \binom{k}{i} \mathrm{Ad}^k ( (1+tA)^{-m_1} ) \mathrFm{Ad}^{k-i} ( tA ) (1+tA)^{-m_2}
% \\
% =& t \sum_{i = 0}^{k-1} \sum_{m_1+m_2 = M+1} 
%   && \big( -2 k \mathrm{Ad}^{k-1} ((1+tA)^{-m_1}) \partial_{x_j} (1+tA)^{-m_2} 
%   \\
%   & &&+ k(k-1) \mathrm{Ad}^{k-2} ((1+tA)^{-m_1}) (1+tA)^{-m_2} \big).
% \end{aligned}
% \]

The above inequalities show that 
\begin{equation}\label{eq:0131-2}
\begin{split}
\| \nabla ^\ell \varphi(tA) f \|_{\ell^1 (L^2)_t}
\leq&\, C(t^{-\frac{\ell}{2}} + t^{-\frac{d}{4}} (t ^{\frac{\alpha-\ell}{2}})^{\frac{d}{2\alpha}}
(t^{-\frac{\ell}{2}})^{1-\frac{d}{2\alpha}}) 
\| (1+tA)^{-\beta} f \|_{\ell^1 (L^2)_t} 
\\
=&\, C t^{-\frac{\ell}{2}} 
\| (1+tA)^{-\beta} f \|_{\ell^1 (L^2)_t} .
\end{split}
\end{equation}
Finally, if we take $\beta$ to be greater than $d/4$, we have a resolvent estimate  from Proposition~4.1 in~\cite{paper:IwMaTa-2018}. 
\[
\| (1 + tA)^{-\beta} f \|_{\ell^1 (L^2)_t}
\leq C t^{-\frac{d}{4}} \| f \|_{L^1}.  
\]
We then obtain by \eqref{eq:0131-1} that 
\[
t^{\frac{d}{4}} \| \nabla ^\ell \varphi(tA) f \|_{\ell^1 (L^2)_t}
\leq C t^{\frac{d}{4}} t^{-\frac{\ell}{2}} \| (1+tA)^{-\beta} f \|_{\ell^1 (L^2)_t}
\leq C t^{-\frac{\ell}{2}} \| f \|_{L^1}. 
\]

% Section3 L^infty estimate
\section{$L^\infty$ estimate in Theorem~\ref{thm:p}}

We prove the estimate for $p = \infty$ in Theorem~\ref{thm:p}, we will show the following.

\begin{thm}
    Let $\ell \geq 2$. Suppose that $\Omega \subset {\mathbb R}^d$ be a bounded domain with $C^{\ell+2}$ boundary.
    Then, there exists $C > 0$ and $T > 0$ such that
    \[
    \sup_{0 < t < T}{\| N_\ell[u](t) \|_{L^\infty}} < C \| u_0 \|_{L^\infty}
    \]
    holds for any solution $u(t)$ of $\eqref{eq:main}$ with $u_0 \in C_c^{\infty}(\Omega)$, where $ N_\ell[u]$ is defined by~\eqref{def:20241216-1}.
\end{thm}

\begin{pr}
    We prove by contradiction. Assume that there exist a sequence $\{ u_m \}_{m \in {\mathbb N}}$ of solution to $\eqref{eq:main}$ and a sequence $\{ \tau_m \}_{m \in {\mathbb N}}$ such that
    \[
    \begin{aligned}
        \| N_\ell[u](\tau_m) \|_{L^\infty} 
        > m \| u_{0, m} \|_{L^\infty}, 
        \\
        \tau_m \searrow 0 \, (m \to \infty),
    \end{aligned}
    \]
    where $u_{0, m} \in C_c^{\infty}(\Omega)$ is an initial data of $u_m$.

    Now we introduce a function $v_m$ as follows.
    It holds that
    $\displaystyle M_m := \sup_{t \in (0, \tau_m)}{\| N[u_m](t) \|_{L^\infty}} < \infty$ 
    by Lemma $\ref{lemma:C2}$, and we define $\widetilde{u}_m := u_m / M_m$.
    Then, there exists some point $(t_m, x_m) \in (0, \tau_m) \times \Omega$ satisfying
    $N_\ell[\widetilde{u}_m](t_m, x_m) \geq 1/2$.
    We therefore define $v_m$ by translating and scaling where the point $(t_m, x_m)$ corresponds to the point $(1, 0)$.
    \begin{align*}
        v_m (t, x) := \widetilde{u}_m (t_m t, x_m + t_m^{1/2} x) 
            \quad \text{ for } (t, x) \in (0, 1) \times \Omega_m, \\
        \text{where } \Omega_m := \Big\{ y \in \mathbb{R}^d \, \Big| \, y = \frac{x - x_m}{t_m^{1/2}}, x \in \Omega \Big\}.
    \end{align*}

    Here $v_m$ satisfies the heat equation and the following inequities.
    \begin{equation} \label{eq:1105-01}
    % \mathrm{N} or N
    \begin{aligned}
        &N_\ell[v_m](1, 0)
        = N_\ell[\widetilde{u}_m](t_m, x_m) \geq \frac{1}{2}, 
        \\
        &\sup_{t \in (0, 1)} \| N_\ell[v_m](t) \|_{L^\infty} 
        = \sup_{t \in (0, t_m)} \| N_\ell[\widetilde{u}_m](t) \|_{L^\infty} \leq 1, 
        \\
        &\| v_{0, m} \|_{L^\infty}  
        = \| \widetilde{u}_{0, m} \|_{L^\infty} < \frac{1}{m}.
    \end{aligned}
    \end{equation}

    The situation depends on whether the limit superior of
    \[
    c_m := \frac{d(\Omega, x_m)}{t_m^{1/2}} = d(\Omega_m, 0)
    \]
    goes infinity or not. In the following, the same indices will be used when taking a subsequence.

    \vspace{3mm}
    
    Case 1. $\displaystyle \limsup_{m \to \infty}c_m = \infty$

    Take a subsequence so that $\displaystyle \lim_{m \to \infty} c_m = \infty$ by the definition of the limit superior.
    % Ascoli-Arzel\`{a} theorem implies that for any compact set $K \subset \mathbb{R}^d$, there exist a subsequence $\{v_{m}\}_{m \in \mathbb{N}}$ and $v^K \in W^{\ell, \infty}(K)$ such that
    % \[
    % v_{m} \to v^K \text{ in } W^{\ell, \infty}(K) \quad (m \to \infty).
    % \]
    % by a compact empbedding $W^{3, r}(K) \subset \subset W^{2, \infty}(K)$.
    For any compact set $K \subset \mathbb{R}^d \times (0, 1]$, $\partial_t v_m = \Delta v_m$ and $\nabla v_m$ are bounded in $K$ by~\eqref{eq:1105-01}. 
    Ascoli-Arzel\`{a} theorem implies that there exist a subsequence $\{v_{m}\}_{m \in \mathbb{N}}$ and $v \in W^{\ell, \infty}(\mathbb{R}^d)$ satisfying
    \begin{equation*}
    v_{m} \to v
    \quad \text{ locally uniformly in } (0, 1] \times \mathbb{R}^d.
    \end{equation*}
    In addition to this, for any compact set $K_X \subset \mathbb{R}^d$ and $t \in (0, 1]$, 
    $| \nabla^{\ell+1} v_m(t) |$ is bounded on $K_X$. Therefore, we also have
    \begin{equation*}
    \left.
        \begin{aligned}
        \nabla v_{m} \to \nabla v&
        \\
        \vdots&
        \\
        \nabla^\ell v_m \to \nabla^\ell v&
        \end{aligned}
    \right \} 
    \quad \text{ locally uniformly in } \mathbb{R}^d \text{ for each } t \in (0, 1].
    \end{equation*}
    
    % Hence, there exist a subsequence $\{v_{m}\}_{m \in \mathbb{N}}$ and $v \in W^{\ell, \infty}(\mathbb{R}^d)$ satisfying
    % \begin{equation*}
    % \begin{aligned}
    %     \left.
    %     v_{m} \to v
    %     \right. &
    %     \quad \text{ locally uniformly in } \mathbb{R}^d \times (0, 1], 
    %     \\
    %     \left.
    %         \begin{aligned}
    %         \nabla v_{m} \to \nabla v&
    %         \\
    %         \vdots&
    %         \\
    %         \nabla^\ell v_m \to \nabla^\ell v&
    %         \end{aligned}
    %     \right \} &
    %     \quad \text{ locally uniformly in } \mathbb{R}^d \text{ for each } t \in (0, 1].
    % \end{aligned}
    % \end{equation*}
    % since for each $\epsilon > 0, r \in (d, \infty)$, 
    % \[
    % \sup_{m\in\mathbb{N}} \sup_{t \in [\epsilon, 1]} \| v_m(t) \|_{W^{\ell+1, r}}
    % \]
    % is bounded by~\eqref{eq:1105-01} and compact embedding $W^{\ell+1, r}(K) \subset\subset W^{\ell, \infty}(K)$ holds.
    
    For any $\phi \in C_c^{\infty} ({\mathbb R}^d \times [0, 1))$, it follows
    \[
        \int_0^1 \int_{{\mathbb R}^d} v_{m} \cdot (\partial_t \phi + \Delta \phi) {\mathrm d}x {\mathrm d}t
        = \int_{{\mathbb R}^d} v_{m} (x, 0) \cdot \phi(x, 0) \, {\mathrm d}x
        \text{ for sufficient large } m, 
    \]
    and by taking the limit as $m \to \infty$, we obtain
    \[
    \int_0^1 \int_{{\mathbb R}^d} v \cdot (\partial_t \phi + \Delta \phi) {\mathrm d}x {\mathrm d}t
    = 0
    \]
    by~\eqref{eq:1105-01} and Lebesgue's dominated convergence theorem.
    Therefore the limit $v$ is the weak solution of heat equation and 
    it follows $v \equiv 0$ by the uniqueness of the heat equation (see section~4.4.2 in~\cite{book:Giga-Giga-Saal}).
    However, it contradicts $\displaystyle N_\ell[v](1, 0) = \lim_{m \to \infty} N_\ell[v_m](1, 0) \geq 1/2$.

    \vspace{3mm}

    Case 2. $\displaystyle \limsup_{m \to \infty}c_m < \infty$

    By taking a subsequence, we assume that $\displaystyle \lim_{m \to \infty} c_m = c_0 < \infty$ and 
    $x_m$ converges to $\hat{x} \in \partial\Omega$.
    In addition, by rotation and translation, we also assume $\hat{x} = 0$.

    Since $\Omega$ is of $C^{\ell+1}$, there exists a $C^{\ell+1}$ function $h$ such that 
    % \[
    % \begin{aligned}
    %     (\Omega_m)_{\text{loc}} := 
    %     \{ (y', y_d) \in \mathbb{R}^d \mid
    %     &h(x_m' + t_m^{1/2} y') < (x_m)_d + t_m^{1/2} y_d
    %     \\
    %     &< h(x_m' + t_m^{1/2} y') + \beta, 
    %     |t_m^{1/2} y'| < \alpha \}
    %     \subset \Omega_m, 
    %     \\
    %     \nabla{h}(0) = 0, h(0) = 0. 
    % \end{aligned}
    % \]
    \[
    \begin{aligned}
        (\Omega_m)_{\text{loc}} := 
        \{ (y', y_d) \in \mathbb{R}^d \mid
        &h_m(y') < y_d
        < h_m(y') + \beta / t_m^{1/2}, 
        |y'| < \alpha / t_m^{1/2} \}
        \subset \Omega_m, 
    \end{aligned}
    \]
    \[
    \displaystyle 
    \begin{aligned}
        \text{ with } h_m(y') = \frac{ h(x_m' + t_m^{1/2} y') - (x_m)_d }{ t_m^{1/2} }, 
        \nabla'{h}(\hat{x}) = 0, h(\hat{x}) = 0, 
    \end{aligned}
    \]
    where $\alpha, \beta$ are the constants independent of $\hat{x} \in \partial\Omega$. 

    Since $d(\Omega_m, x_m) / (x_m)_d \to 1$ as $m \to \infty$, 
    the domain $(\Omega_m)_{\text{loc}}$ approaches to 
    \[
    \mathbb{R}^d_{+, -c_0} := \{ (x', x_d) \in \mathbb{R}^d \mid x_d > -c_0 \}.
    \]
    % We check that $h_m(\cdot)$ converges to $-c_0$ locally uniformly up to $(\ell+1)$-th derivatives: 
    % \[
    % \begin{aligned}
    % \sup_{|y'| < r} \partial_{y_i} h_m(y') 
    % = \sup_{|y'| < r} \partial_{y_i} h(x_m' + t_m^{1/2} y')
    % \to \partial_{y_i} h(\hat{x}) = 0, 
    % \\
    % \sup_{|y'| < r} \partial_{y'}^{\gamma} h_m(y') 
    % = t_m^{\frac{|\gamma| - 1}{2}} \sup_{|y'| < r} \partial_{y'}^{\gamma} h(x_m' + t_m^{1/2} y')
    % \to 0
    % \end{aligned}
    % \]
    % as $m \to \infty$, for any $r > 0$ and multi-index $\gamma$ with $2 \leq |\gamma| \leq \ell+1$.

    By the similar argument as in Case 1, there exists a subsequence $\{v_{m}\}_{m \in \mathbb{N}}$ and $v \in W^{\ell, \infty}(\mathbb{R}^d)$ such that
    \begin{equation*}
    \begin{aligned}
        \left.
        v_{m} \to v
        \right. &
        \quad \text{ locally uniformly in } \overline{\mathbb{R}}^d_{+, -c_0} \times (0, 1], 
        \\
        \left.
            \begin{aligned}
            \nabla v_{m} \to \nabla v&
            \\
            \vdots&
            \\
            \nabla^\ell v_m \to \nabla^\ell v&
            \end{aligned}
        \right \} &
        \quad \text{ locally uniformly in } \overline{\mathbb{R}}^d_{+, -c_0} \text{ for each } t \in (0, 1].
    \end{aligned}
    \end{equation*}

    Fix $R > 0$ and let $B_R^+ = B_R(0, \ldots, 0, -c_0) \cap {\mathbb R}^d_{+, -c_0}$.
    Set a test funciton 
    \[
    \phi \in  C\Big([0, 1] ; W^{1,1}_0(R^d_{+, -c_0}) \cap W^{2, 1}(R^d_{+, -c_0}) \cap C(R^d_{+, -c_0})\Big) \cap C^1\Big([0, 1] ; L^1(R^d_{+, -c_0})\Big)
    \]
    with $|\nabla \phi(x)| \leq C (1 + |x|^2)^{-d/2}$, then $v_m$ satisfies
    \[
    0 
    = \int_0^1 \int_{\Omega_m \cap B_R^+} \Delta {v}_m \phi \, {\mathrm d}x {\mathrm d}t
    - \int_0^1 \int_{\Omega_m \cap B_R^+} (\partial_t {v}_m) \phi \, {\mathrm d}x {\mathrm d}t.
    \]

    By following the proof of Theorem~4.1 in~\cite{paper:FuIwKo-2024}, we then obtain
    \[
    \begin{aligned}
    &\lim_{R \to \infty} \lim_{m \to \infty} \int_0^1 \int_{\Omega_m \cap B_R^+} \Delta {v}_m \phi \, {\mathrm d}x {\mathrm d}t
    = \int_0^1 \int_{\mathbb R^d_{+, -c_0}} {v} \Delta \phi \, {\mathrm d}x {\mathrm d}t, 
    \\
    &\lim_{R \to \infty} \lim_{m \to \infty} \int_0^1 \int_{\Omega_m \cap B_R^+} \partial_t {v}_m \phi \, {\mathrm d}x {\mathrm d}t
    = - \int_0^1 \int_{\mathbb R^d_{+, -c_0}} {v} \partial_t \phi \, {\mathrm d}x {\mathrm d}t
    - \int_{\mathbb R^d_{+, -c_0}} {v}(0, x) \phi(0, x) \, {\mathrm d}x
    \end{aligned}
    \]
    and therefore, ${v}$ satisfies
    \[
    \int_0^1 \int_{\mathbb R^d_{+, -c_0}} {v} (\Delta \phi +\partial_t \phi) \, {\mathrm d}x {\mathrm d}t
    = - \int_{\mathbb R^d_{+, -c_0}} {v}(0, x) \phi(0, x) \, {\mathrm d}x.
    \]
    
    By the uniqueness of the heat equation for the half space (Proposition~\ref{prop:0126-1}), it follows ${v} \equiv 0$.
    However, as in Case 1, it contradicts $N[{v}](1, 0) \geq 1/2$.
\end{pr}

% Theorem 1.2
\section{Proof of Theorem~\ref{thm:20250127-1}}
Let $\{ \phi_j \}_{j \in \mathbb{Z}}$ be a function sequence satisfying the following property and consider dyadic decomposition.
\begin{equation*}
\begin{aligned}
\phi_j \in C_c^\infty((0, \infty)), \quad
\text{supp }\phi_j \subset [2^{j-1}, 2^{j+1}], 
\\
\phi_j(\cdot) = \phi_0(2^{-j} \cdot), \quad
\sum_{j \in \mathbb{Z}} \phi_j(\xi) = 1, \, \xi \in (0, \infty).
\end{aligned}
\end{equation*}
Let $f \in C_c^\infty(\Omega)$ .
For each $j \in \mathbb{Z}$, define $\phi_j(\sqrt{A})$ as an operator on $L^2(\Omega)$ such that
\begin{equation*}
\phi_j(\sqrt{A}) = \int_{0}^{\infty} \phi_j ( \sqrt{\lambda} ) \, \mathrm{d}E_A(\lambda).
\end{equation*}
Then, note that
\begin{equation*}
\begin{aligned}
\displaystyle
f
=&\, \sum_{j \in \mathbb{Z}} \phi_j( \sqrt{A} ) f
% \\
% =&\, \sum_{j \in \mathbb{Z}} \left( \phi_{j-1}(\sqrt{A}) + \phi_{j}(\sqrt{A}) + \phi_{j+1}(\sqrt{A}) \right) \phi_j(\sqrt{A}) f 
\quad \text{ in } L^2(\Omega).
\end{aligned}
\end{equation*}

Let $k \in \mathbb{Z}_{\geq 0}$ be a derivative order.
For the $L^p$ norm of the $k$-th order of $e^{-tA^{\alpha / 2}} f$, it holds that
\begin{equation} \label{eq:20250117-4}
\begin{aligned}
\displaystyle
\| \nabla^k e^{-tA^{\alpha / 2}} f \|_{L^p(\Omega)}
=&\, \left\| \sum_{j \in \mathbb{Z}} \nabla^k e^{-tA^{\alpha / 2}} \phi_j(\sqrt{A}) f \right\|_{L^p(\Omega)}
\\
\leq&\, \sum_{j \in \mathbb{Z}} \| \nabla^k e^{-2^{-2j}A} e^{2^{-2j}A} e^{-tA^{\alpha / 2}} \phi_j(\sqrt{A}) f \|_{L^p(\Omega)}
\\
\leq&\, \sum_{j \in \mathbb{Z}} \| \nabla^k e^{-2^{-2j}A} \|_{L^p \to L^p}
\| e^{2^{-2j}A} e^{-tA^{\alpha / 2}} \phi_j(\sqrt{A}) \|_{L^p \to L^p} \| f \|_{L^p(\Omega)}.
\end{aligned}
\end{equation}
Theorem~\ref{thm:p} implies 
\[
\| \nabla^k e^{-2^{-2j}A} \|_{L^p \to L^p} \leq C ( 2^{-2j} )^{-k/2} = C 2^{jk}.
\]
Also, since the spectrum is restricted to a binary number for each $j$, we apply Lemma~2.1 in~\cite{paper:Iw-2018} and then
\begin{equation*}
\begin{aligned}
\| e^{2^{-2j}A} e^{-tA^{\alpha / 2}} \phi_j(\sqrt{A}) \|_{L^p \to L^p}
\leq&\, \| e^{\lambda} e^{-t2^{\alpha j} \lambda^{\alpha / 2}} \phi_0(\sqrt{\lambda}) \|_{H^s(\mathbb{R})}
\\
\leq&\, C e^{-ct 2^{\alpha j}}, 
\end{aligned}
\end{equation*}
where $s$ is a real number such that $s > (d+1) / 2$. 
Therefore, we estimate \eqref{eq:20250117-4} as follows.
\begin{equation*}
\begin{aligned}
\displaystyle
\| \nabla^k e^{-tA^{\alpha / 2}} f \|_{L^p(\Omega)}
\leq&\, \sum_{j \in \mathbb{Z}} C 2^{jk} e^{-ct 2^{\alpha j}} \| f \|_{L^p(\Omega)}
\\
=&\, C t^{-\frac{k}{\alpha}} \sum_{j \in \mathbb{Z}} ( t^{\frac{1}{\alpha}} 2^{j})^k e^{-ct 2^{\alpha j}} \| f \|_{L^p(\Omega)}.
\end{aligned}
\end{equation*}
Here, consider the boundedness of the following quantity.
\[
I(t) \coloneqq \sum_{j \in \mathbb{Z}} ( t^{\frac{1}{\alpha}} 2^{j})^k e^{-ct 2^{\alpha j}}, \quad t > 0.
\]
For any $t > 0$, it follows 
\begin{equation*}
\begin{aligned}
I(2^\alpha t)
=&\, \sum_{j \in \mathbb{Z}} ( (2^\alpha t)^{\frac{1}{\alpha}} 2^{j})^k e^{- c(2^\alpha t) 2^{\alpha j}}
\\
=&\, \sum_{j \in \mathbb{Z}} ( t^{\frac{1}{\alpha}} 2^{j+1})^k e^{- ct 2^{\alpha (j+1)}}
\\
=&\, I(t).
\end{aligned}
\end{equation*}
Furthermore, since $I(t)$ is continuous with respect to $t$, it satisfies
\begin{equation*}
\begin{aligned}
I(t) \leq \max_{1 \leq s \leq 2^\alpha} I(s) < \infty, \quad t > 0.
\end{aligned}
\end{equation*}

Summarising the above, we obtain
\begin{equation*}
\begin{aligned}
\displaystyle
\| \nabla^k e^{-tA^{\alpha / 2}} f \|_{L^p(\Omega)}
\leq\, C t^{-\frac{k}{\alpha}} I(t) \| f \|_{L^p(\Omega)}
\leq\, C t^{-\frac{k}{\alpha}} \| f \|_{L^p(\Omega)}, 
\end{aligned}
\end{equation*}
which proves Theorem~\ref{thm:20250127-1}.

\appendix

\section{Estimate of commutator for resolvent} \label{sec:20241217-4}

We show the following lemma used in the proof for $p = 1$. 
Define
\[
R_{M, t} := (1 + tA)^{-M}
\]
and the sequence $\{ \mathrm{Ad}^k(L) \}_{k \geq 0}$ of operators as follows.
\[
\begin{aligned}
\mathrm{Ad}^0 (L) &:=  L, 
\\
\mathrm{Ad}^k (L) &:= \mathrm{Ad}^{k-1} \left( (x_j - t^{1/2}n_j) L - L (x_j - t^{1/2}n_j) \right),  
\end{aligned}
\]
where $L$ is a operator on $L^2(\Omega)$, and $x_j$ and $n_j$ denote the $j$-th component of $x$ and $n$.

The following lemma for $M = 1$ is known in the paper by Jensen-Nakamura~\cite{paper:JeNa-1995}.

\begin{lem} \label{lem:20241217-2}
Let $M \geq 0$ and $\ell > 0$ be integers and $t > 0$. 
Then, it holds 
\begin{equation} \label{eq:20241210-7}
\mathrm{Ad}^\ell ( R_{M, t} )
= \sum_{k = 0}^{\ell-1} \sum_{m_1 + m_2 = M + 1}
\binom{\ell}{k} \mathrm{Ad}^k ( R_{m_1, t} ) \mathrm{Ad}^{\ell-k} ( tA ) R_{m_2, t}.
\end{equation}
where $m_1$ and $m_2$ in the above sum satisfies $m_1, m_2 \geq 1$ and $\displaystyle \binom{\ell}{k}$ is a binomal coefficient
\[
\begin{aligned}
\displaystyle
\binom{\ell}{k} = 
\begin{cases}
\dfrac{\ell !}{k! (\ell - k)!} & \textup{ if } 0 \leq k \leq \ell
\\
0 & \textup{ otherwise } .
\end{cases}
\end{aligned}
\]
\end{lem}
\begin{pr}
We prove by induction. When $\ell =  1$, we have
\[
\begin{aligned}
\mathrm{Ad}^1 ( R_{M, t} ) 
=&\, (1+tA)^{-M}(x_j - t^{1/2} n_j) - (x_j - t^{1/2} n_j) (1+tA)^{-M}
\\
=&\, \sum_{i = 0}^{M-1} \Big( (1+tA)^{-(M-i)} (x_j - t^{1/2} n_j) (1+tA)^{-i}
  \\
  &\quad\quad\quad - (1+tA)^{-(M-i-1)} (x_j - t^{1/2} n_j) (1+tA)^{-(i+1)} \Big)
\\
=&\, \sum_{i = 0}^{M-1} R_{M-i, t} \left[ (x_j - t^{1/2} n_j) (1 + tA) -  (1 + tA) (x_j - t^{1/2} n_j) \right] R_{i+1, t}
\\
=&\, \sum_{i = 0}^{M-1} \mathrm{Ad}^0 ( R_{M-i, t} ) \mathrm{Ad}^{1} ( tA ) R_{i+1, t}
\\
=&\, \sum_{m_1+m_2 = M+1} \mathrm{Ad}^0 ( R_{m_1, t} ) \mathrm{Ad}^{1} ( tA ) R_{m_2, t}.
\end{aligned}
\]
Let us suppose that~\eqref{eq:20241210-7} is true for $\ell \in \{ 1, 2, \ldots, s \}$. 
Then, for $\ell = s+1$, 
\begin{equation} \label{eq:20250115-1}
\begin{aligned}
\mathrm{Ad}^{s+1} ( R_{M, t} ) 
=&\, \mathrm{Ad}^{1} ( \mathrm{Ad}^{s} ( R_{M, t} ) )
\\
=&\, \mathrm{Ad}^{1} \left( \sum_{k = 0}^{s-1} \sum_{m_1 + m_2 = M + 1}
\binom{s}{k} \mathrm{Ad}^k ( R_{m_1, t} ) \mathrm{Ad}^{s-k} ( tA ) R_{m_2, t} \right)
\\
=&\, \sum_{k = 0}^{s-1} \sum_{m_1+m_2 = M+1} 
  \binom{s}{k} \Big( \mathrm{Ad}^{k+1} ( R_{m_1, t} ) \mathrm{Ad}^{s-k} ( tA ) R_{m_2, t}
  \\
  & \hspace{40mm} + \mathrm{Ad}^{k} ( R_{m_1, t} ) \mathrm{Ad}^{s-k+1} ( tA ) R_{m_2, t}
  \\
  & \hspace{40mm} + \mathrm{Ad}^{k} ( R_{m_1, t} ) \mathrm{Ad}^{s-k} ( tA ) \mathrm{Ad}^1(R_{m_2, t}) \Big)
\\
\eqqcolon& \, I_1 + I_2 + I_3.
\end{aligned}
\end{equation}
For $I_1$, converting $k+1 = k'$,  we write
\begin{equation*}
\begin{aligned}
I_1
=& \, \sum_{m_1+m_2 = M+1} \sum_{k' = 1}^{s}  
\binom{s}{k'-1} \Big( \mathrm{Ad}^{k'} ( R_{m_1, t} ) \mathrm{Ad}^{s-k'+1} ( tA ) R_{m_2, t} \Big)
\\
=& \, \sum_{m_1+m_2 = M+1} \Bigg[ \sum_{k' = 0}^{s} 
\binom{s}{k'-1} \Big( \mathrm{Ad}^{k'} ( R_{m_1, t} ) \mathrm{Ad}^{s-k'+1} ( tA ) R_{m_2, t} \Big)
  \\
  &\, \hspace{25mm}  - \binom{s}{-1} \mathrm{Ad}^{0} ( R_{m_1, t} ) \mathrm{Ad}^{s+1} ( tA ) R_{m_2, t} \Bigg]
\\
=& \, \sum_{m_1+m_2 = M+1} \sum_{k' = 0}^{s} 
\binom{s}{k'-1} \Big( \mathrm{Ad}^{k'} ( R_{m_1, t} ) \mathrm{Ad}^{s-k'+1} ( tA ) R_{m_2, t} \Big).
\end{aligned}
\end{equation*}
Also, for $I_2$, by adding and subtracting the terms $k = s$, we obtain 
\begin{equation*}
\begin{aligned}
I_2
=& \, \sum_{m_1, m_2} \Bigg[ \sum_{k = 0}^{s} \binom{s}{k} \mathrm{Ad}^{k} ( R_{m_1, t} ) \mathrm{Ad}^{s-k+1} ( tA ) R_{m_2, t} 
  \\
  &\, \hspace{15mm} - \binom{s}{s} \mathrm{Ad}^s(R_{m_1, t}) \mathrm{Ad}^{1} ( tA ) R_{m_2, t} \Bigg].
\end{aligned}
\end{equation*}
Therefore, 
\begin{equation*}
\begin{aligned}
&I_1 + I_2 + I_3
\\
=& \sum_{m_1, m_2} \sum_{k = 0}^{s} \binom{s}{k-1} \mathrm{Ad}^{k} ( R_{m_1, t} ) \mathrm{Ad}^{s-k+1} ( tA ) R_{m_2, t}  
  \\
  &\, + \sum_{m_1, m_2} \left[ \sum_{k = 0}^{s} \binom{s}{k} \mathrm{Ad}^{k} ( R_{m_1, t} ) \mathrm{Ad}^{s-k+1} ( tA ) R_{m_2, t} - \mathrm{Ad}^s(R_{m_1, t}) \mathrm{Ad}^{1} ( tA ) R_{m_2, t} \right]
  \\
  &\, + \sum_{m_1, m_2} \sum_{k = 0}^{s-1} \binom{s}{k} \mathrm{Ad}^{k} ( R_{m_1, t} ) \mathrm{Ad}^{s-k} ( tA ) \mathrm{Ad}^1(R_{m_2, t}).
\end{aligned}
\end{equation*}
Here, using the relation
\[
\binom{s}{k-1} + \binom{s}{k} = \binom{s+1}{k}
, \quad k, s \in \mathbb{Z}_{\geq 0}
\]
for the binomial coefficients, the right-hand side of~\eqref{eq:20250115-1} can be expressed as following.
\begin{equation*}
\begin{aligned}
\mathrm{Ad}^{s+1} ( R_{M, t} ) 
=& \, I_1 + I_2 + I_3
\\
=& \sum_{m_1, m_2} \sum_{k = 0}^{s} \binom{s+1}{k} \mathrm{Ad}^{k} ( R_{m_1, t} ) \mathrm{Ad}^{s-k+1} ( tA ) R_{m_2, t}
  \\
  &\, - \sum_{m_1, m_2}  \mathrm{Ad}^s(R_{m_1, t}) \mathrm{Ad}^{1} ( tA ) R_{m_2, t}
  \\
  &\, + \sum_{m_1, m_2} \sum_{k = 0}^{s-1} \binom{s}{k} \mathrm{Ad}^{k} ( R_{m_1, t} ) \mathrm{Ad}^{s-k} ( tA ) \mathrm{Ad}^1(R_{m_2, t})
\\
=:& \, I + II + III.
\end{aligned}
\end{equation*}
By the induction assumption, $II$ and $III$ can be written as
\[
\begin{aligned}
II =& - \sum_{m_1, m_2} \left( \sum_{k=0}^{s-1} \sum_{m_3 + m_4 = m_1 + 1} 
\binom{s}{k} \mathrm{Ad}^k ( R_{m_3, t} ) \mathrm{Ad}^{s-k} ( tA ) R_{m_4, t} 
\right) \mathrm{Ad}^{1} ( tA ) R_{m_2, t}
\\
=& - \sum_{k = 0}^{s-1} \sum_{m_3 + m_4 + m_2 = M+2} 
\binom{s}{k} \mathrm{Ad}^k ( R_{m_3, t} ) \mathrm{Ad}^{s-k} ( tA ) R_{m_4, t} 
 \mathrm{Ad}^{1} ( tA ) R_{m_2, t}, 
\end{aligned}
\]
\[
\begin{aligned}
III =& \sum_{m_1, m_2} \sum_{k = 0}^{s-1} \binom{s}{k} \mathrm{Ad}^{k} ( R_{m_1, t} ) \mathrm{Ad}^{s-k} ( tA ) 
\left( \sum_{m_5 + m_6 = m_2 + 1} R_{m_5, t} \mathrm{Ad}^{1} ( tA ) R_{m_6, t} \right)
\\
=& \sum_{k = 0}^{s-1} \sum_{m_1 + m_5 + m_6 = M+2} 
\binom{s}{k} \mathrm{Ad}^k ( R_{m_1, t} ) \mathrm{Ad}^{s-k} ( tA ) R_{m_5, t} 
\mathrm{Ad}^{1} ( tA ) R_{m_6, t}.
\end{aligned}
\]
Therefore we have $II + III = 0$ and
\[
\begin{aligned}
\mathrm{Ad}^{s+1} ( R_{M, t} ) 
=&\, I 
\\
=&\, \sum_{m_1, m_2} \sum_{k = 0}^{s} \binom{s+1}{k} \mathrm{Ad}^{k} ( R_{m_1, t} ) \mathrm{Ad}^{s+1-k} ( tA ) R_{m_2, t}, 
\end{aligned}
\]
which implies \eqref{eq:20241210-7} for $\ell = s + 1$.
\end{pr}

\begin{lem} \label{lem:20241217-3}
Let $M \geq 0$. Then, for any $t, \tau \in \mathbb{R}$, 
\begin{equation}
	\mathrm{Ad}^1(e^{-i \tau R_{M, t}}) 
    = - i \int^{\tau}_{0} e^{-isR_{M, t}} \mathrm{Ad}^1(R_{M, t}) e^{-i(\tau-s)R_{M, t}} 
      \,{\mathrm d}s.
\end{equation}
Moreover, for $k \geq 0$, the following holds.
\begin{equation} \label{eq:20241217-1}
\begin{aligned}
	&\mathrm{Ad}^{k+1}(e^{-itR_{M, t}}) 
    \\
	=& - i \int^{\tau}_{0} \sum_{k_1 + k_2 + k_3 = k}
	\frac{k!}{k_1! k_2! k_3!} \mathrm{Ad}^{k_1}(e^{-isR_{M, t}}) \mathrm{Ad}^{k_2+1}(R_{M, t}) \mathrm{Ad}^{k_3}(e^{-i(\tau-s)R_{M, t}}) \, {\mathrm d}s, 
\end{aligned}
\end{equation}
where $k_1, k_2, k_3$ in the above sum satisfies $k_1, k_2, k_3 \geq 0$.
\end{lem}

\begin{pr}
We prove by induction with respect to $k$. When $k = 0$, we have
\[
\begin{aligned}
\mathrm{Ad}^1(e^{-i \tau R_{M, t}})
& = x_j e^{-i \tau R_{M, t}} - e^{-i \tau R_{M, t}} x_j 
\\
& = - \int^\tau_0 \frac{d}{ds}
\big(e^{-isR_{M, t}} \, x_j \,  e^{-i(\tau-s)R_{M, t}} \big)\, {\mathrm d}s
\\
& =
-i \int^\tau_0 e^{-isR_{M, t}}
(x_j R_{M, t} - R_{M, t} x_j) e^{-i(\tau-s)R_{M, t}} \, {\mathrm d}s 
\\
& =
-i \int^\tau_0 e^{-isR_{M, t}}
\mathrm{Ad}^1 (R_{M, t}) e^{-i(\tau-s)R_{M, t}} \, {\mathrm d}s.
\end{aligned}
\]

Let us suppose that~\eqref{eq:20241217-1} is true for $k \in \{ 0, 1, \ldots, l \}$. 
Then, for $k = l+1$, 
\begin{align*}
& \mathrm{Ad}^{l+2}(e^{-i \tau R_{M, t}}) 
\\
=& (x_j - t^{1/2}n_j) \mathrm{Ad}^{l+1}(e^{-i \tau R_{M, t}}) 
  - \mathrm{Ad}^{l+1}(e^{-i \tau R_{M, t}}) (x_j - t^{1/2}n_j) 
\\
= &
    - i \int^{\tau}_{0} \sum_{k_1 + k_2 + k_3 = l}
    \frac{l!}{k_1! k_2! k_3!} \mathrm{Ad}^{k_1+1}(e^{-isR_{M, t}}) \mathrm{Ad}^{k_2+1}(R_{M, t}) \mathrm{Ad}^{k_3}(e^{-i(\tau-s)R_{M, t}}) \, {\mathrm d}s
    \\
    &
    - i \int^{\tau}_{0} \sum_{k_1 + k_2 + k_3 = l}
    \frac{l!}{k_1! k_2! k_3!} \mathrm{Ad}^{k_1}(e^{-isR_{M, t}}) \mathrm{Ad}^{k_2+2}(R_{M, t})  \mathrm{Ad}^{k_3}(e^{-i(\tau-s)R_{M, t}}) \, {\mathrm d}s
    \\
    &
    - i \int^{\tau}_{0} \sum_{k_1 + k_2 + k_3 = l}
    \frac{l!}{k_1! k_2! k_3!} \mathrm{Ad}^{k_1}(e^{-isR_{M, t}}) \mathrm{Ad}^{k_2+1}(R_{M, t}) \mathrm{Ad}^{k_3+1}(e^{-i(\tau-s)R_{M, t}}) \, {\mathrm d}s
\\
= &
    - i \int^{\tau}_{0} \sum_{k_1 + k_2 + k_3 = l+1}
    \frac{k_1}{l+1} \frac{(l+1)!}{k_1! k_2! k_3!} 
    \mathrm{Ad}^{k_1}(e^{-isR_{M, t}})
    \mathrm{Ad}^{k_2+1}(R_{M, t})
    \mathrm{Ad}^{k_3}(e^{-i(\tau-s)R_{M, t}}) \, {\mathrm d}s 
    \\
    &
    - i \int^{\tau}_{0} \sum_{k_1 + k_2 + k_3 = l+1}
    \frac{k_2}{l+1} \frac{(l+1)!}{k_1! k_2! k_3!}
    \mathrm{Ad}^{k_1}(e^{-isR_{M, t}})
    \mathrm{Ad}^{k_2+1}(R_{M, t}) 
    \mathrm{Ad}^{k_3}(e^{-i(\tau-s)R_{M, t}}) \, {\mathrm d}s
    \\
    &
    - i \int^{\tau}_{0} \sum_{k_1 + k_2 + k_3 = l+1}
    \frac{k_3}{l+1} \frac{(l+1)!}{k_1! k_2! k_3!} 
    \mathrm{Ad}^{k_1}(e^{-isR_{M, t}}) 
    \mathrm{Ad}^{k_2+1}(R_{M, t}) 
    \mathrm{Ad}^{k_3}(e^{-i(\tau-s)R_{M, t}}) \, {\mathrm d}s
\\
= &
- i \int^{\tau}_{0} \sum_{k_1 + k_2 + k_3 = l+1}
\frac{(l+1)!}{k_1! k_2! k_3!}
\mathrm{Ad}^{k_1}(e^{-isR_{M, t}}) 
\mathrm{Ad}^{k_2+1}(R_{M, t}) 
\mathrm{Ad}^{k_3}(e^{-i(\tau-s)R_{M, t}}) \, {\mathrm d}s .
\end{align*}
Therefore, we get~\eqref{eq:20241217-1} for $k = l+1$ and it proves the lemma.
\end{pr}

\end{document}